\documentclass[12pt]{amsart}
\usepackage{epsfig}
\usepackage{amsmath}
\usepackage{amsfonts}
\usepackage{amssymb}
\usepackage{epigraph, fancyhdr}

\def\D{\mathcal D}

\def\K{\mathcal K}
\def\J{\mathcal J}
\def\L{\mathcal L} 

\def\O{\mathcal O}

\def\1{\mathbf 1}

\def\QQ{\mathbb Q}

\def\CC{\mathbb C}

\def\Res{\operatorname{Res}}

\def\hat{\widehat}

\def\p{\partial}

\def\t{{\mathbf t}}

\def\la{\lambda}
\def\gL{\Lambda}

\def\ev{\operatorname{ev}}
\def\ft{\operatorname{ft}}

\def\Eu{\operatorname{Euler}}
\def\fake{\operatorname{fake}}
\def\tr{\operatorname{tr}}

\renewcommand{\Delta}{\triangle}

\title[$\D_q$-modules]
      {Permutation-equivariant \\ quantum K-theory IV. \\
      $\D_q$-modules}

\author[A. Givental]{Alexander GIVENTAL}
\thanks{This material is based upon work supported by the National 
Science Foundation under Grant DMS-1007164, and by the IBS Center for Geometry 
and Physics, POSTECH, Korea.} 

\date{July 15, 2015}

\begin{document}

\begin{abstract} 
In Part II, we saw how permutation-equivariant quantum K-theory of a manifold with isolated fixed points of a torus action can be reduced via fixed point localization to permutation-equivariant quantum K-theory of the point.
In Part III, we gave a complete description of permutation-equivariant quantum K-theory of the point by means of adelic characterization. Here we apply the adelic characterization to introduce the action on this theory of a certain group of $q$-difference operators. This action enables us to prove that toric $q$-hypergeometric functions represent K-theoretic GW-invariants of toric manifolds. 
\end{abstract}

\maketitle

\section*{Overruled cones and $\D_q$-modules}

In Part III, we gave the following {\em adelic characterization} of the big J-function $\J_{pt}$ of the point target space. In the space $\K$ of ``rational functions'' of $q$ (consisting in fact of series in auxiliary variables with coefficients which are rational functions of $q$), let $\L$ denote the range of $\J_{pt}$. We showed that {\em an element $f\in \K$ lies in $\L$ if and only if Laurent series expansions $f_{(\zeta)}$ of $f$ near $q=\zeta^{-1}$ satisfy 

  (i) $f_{(1)} = (1-q) e^{\tau/(1-q)} \times (\text{\em power series in $q-1$})$ for some $\tau \in \gL_{+}$,\footnote{For convergence purposes, we assume that the Adams operations $\Psi^k$  on $\gL$ with $k>1$ increase certain filtration $\gL\supset \gL_{+}\supset \gL_{++}\supset \cdots $, and that the domain of the J-function is $\gL_{+}$.}

(ii) when $\zeta\neq 1$ is a primitive $m$-th root of unity,
  \[ f_{(\zeta)}(q^{1/m}/\zeta) = \Psi^m(f_{(1)}/(1-q)) \times (\text{\em power series in $q-1$}),\]
where $\Psi^m$ is the Adams operation extended from $\gL$ by $\Psi^m(q)=q^m$;

(iii) when $\zeta \neq 0,\infty$ is not a root of unity, $f_{(\zeta)}(q/\zeta)$ is a power series in $q-1$.}

\medskip

Another way to phrase (i) is to say that $f_{(1)}$ lies in the range $\L^{fake}$ of the {\em ordinary} (or {\em fake}) J-function $\J_{pt}^{ord}$ in the space $\hat{\K}=\gL ((q-1))$ of Laurent series in $q-1$:
\[ \L^{fake}= \bigcup_{\tau \in\gL_{+}}(1-q) e^{\tau/(1-q)} \hat{K}_{+}, \ \ \ \hat{K}_{+}:=\gL[[q-1]].\]
The range $\L^{\fake}$ is an example of an {\em overruled cone}: Its tangent spaces $T_{\tau}=e^{\tau/(1-q)}\hat{\K}_{+}$ are tangent to $\L^{fake}$ along the subspaces $(1-q) T_{\tau}$ (which sweep $\L^{fake}$ as the parameter $\tau$ varies through $\gL_{+}$.\footnote{In terminology of S. Barannikov \cite{Bar}, this is a variation of semi-infinite Hodge structures: The flags $\cdots \subset (1-q)T_{\tau} \subset T_{\tau}\subset (1-q)^{-1}T_{\tau} \subset \cdots $ vary in compliance with ``Griffiths' transversality condition''.}) As it will be explained shortly, this property leads to the invariance of $\L^{fake}$ to  certain finite-difference operators.

Recall that in permutation-equivariant quantum K-theory, we work over a $\lambda$-algebra, a ring equipped with Adams homomorphisms $\Psi^m$, $m=1,2,\dots$, $\Psi^1=\operatorname{Id}$, $\Psi^m\Psi^l=\Psi^{ml}$.
Let us take $\gL := \gL_0 [[\la, Q]]$ with $\Psi^m(\la)=\la^m$,  $\Psi^m(Q)=Q^m$, where $\gL_0$ is any ground $\lambda$-algebra over $\CC $.

Consider the algebra of finite-difference operators in $Q$. Such an operator is a non-commutative expression $D(Q, 1-q^{Q\p_Q}, q^{\pm 1})$. Clearly, the space
$\hat{\K}_{+}=\gL [[q-1]]$ (as well as $(1-q)\hat{\K}_{+}$) is a $\D_q$-module. Consequently each ruling space $(1-q)T_{\tau}=e^{\tau/(1-q)}(1-q)\hat{K}_{+}$ is a $\D_q$-module too. Indeed, 
\[ q^{Q\p_Q}e^{\tau(Q)/(1-q)}=e^{\tau (Q)/(1-q)}e^{(\tau(qQ)-\tau(Q))/(1-q)},\]
where the second factor lies in $\hat{\K}_{+}$. Moreover, we have

\medskip

{\tt Proposition.} {\em $e^{\la D(Q,1-q^{Q\p_Q},q)/(1-q)} \L^{fake} = \L^{fake}$.}

\medskip



{\tt Proof.} The ruling space $(1-q) T_{\tau}$ is a $\D_q$-module, and hence invariant under $D$. Therefore for $f\in (1-q)T_{\tau}$, we have $Df/(1-q)\in T_{\tau}$, i.e. the vector field defining the flow $t \mapsto e^{t \la D/(1-q)}$ is tangent to $\L^{fake}$, and so the flow preserves $\L^{fake}$. It remains to take $t=1$, which is possible thanks to $\la$-adic convergence.

\medskip

{\tt Remark.} Generally speaking, linear transformation $e^{ \la D/(1-q)}$ does not preserve ruling spaces $(1-q)T_{\tau}$, but transforms each of them into another such space. Indeed, preserving $\L^{fake}$, it 
transform tangent spaces $T_{\tau}$ into tangent spaces, and since it commutes
with multiplication by $1-q$, it also transforms ruling spaces $(1-q)T_{\tau}$
into ruling spaces.

\medskip

Likewise, cone $\L \subset \K$ is ruled by subspaces comparable to $(1-q)\K_{+}$, namely by $(1-q)L_{\tau}$, where $L_{\tau}:=e^{\sum_{k>0}\Psi^k(\tau)/k(1-q^k)}\K_{+}$. However $L_{\tau}$ are not tangent to $\L$. Nonetheless the following result holds. 

\medskip

{\tt Theorem.} {\em The range $\L$ of the big J-function $\J_{pt}$ in the permutation-equivariant quantum K-theory of the point target space is preserved
by operators of the form 
\[ e^{\sum_{k>0} \la^k \Psi^k\left( D(1-q^{kQ\p_Q},q^{\pm 1})\right)/k(1-q^k)}.\]}

\medskip

{\tt Remarks.} (1) The operator $D$ has constant coefficients, i.e. is independent of $Q$.

(2) Note that $\Psi^k(q^{Q\p_Q})=q^{kQ^k\p_{Q^k}}=q^{Q\p_Q}$, and not $q^{kQ\p_Q}$ as in the exponent.

(3) The reader is invited to check that the theoem and its proof are extended without any changes to the case finite difference operators in several variables $Q_1,\dots, Q_K$. We will use the theorem in this more general form in Part V.

\medskip

{\tt Proof.} Assuming that $(1-q)f \in \L$, we use the adelic characterization of $\L$ to show that $(1-q) g \in \L$, where
 \[ g(q) := e^{\sum_{k>0} \la^k \Psi^k\left(D(1-q^{kQ\p_Q},q^{\pm 1})\right)/k(1-q^k)} f(q).\]
First, this relationship between $g$ and $f$ also holds between $g_{(1)}$ and $f_{(1)}$ where however both sides need to be understood as Laurent series in $q-1$. Since $f_{(1)}\in \L^{fake}$, Proposition implies that $g_{(1)} \in \L^{fake}$ too.
 
Next, applying $\Psi^m$ to both sides, we find
\[ \Psi^m(g_{(1)})= e^{\sum_{l>0} \la^{ml} \Psi^{ml}\left(D(1-q^{lQ\p_Q},q^{\pm1})\right)/l(1-q^{ml})} \Psi^m(f_{(1)}).\]
On the other hand, for an $m$-th primitive root of unity $\zeta$, taking into account that $\Psi^{ml}(q)=q^{ml}$ turns after the change $q\mapsto q^{1/m}/\zeta$ into $q^l$, and that $q^{mlQ\p_Q}$ turns after this change into
$q^{lQ\p_Q}$, we find
\[ g_{(\zeta)}(q^{1/m}/\zeta)=e^{\Delta} e^{\sum_{l>0}\la^{ml} \Psi^{ml}\left(D(1-q^{lQ\p_Q},q^{\pm 1/m})\right)/ml(1-q^l)} \, f_{(\zeta)}(q^{1/m}/\zeta),\]
where the finite-difference operator $\Delta$ has coefficients regular at $q=1$. Here we factor off the terms regular at $q=1$ using the fact that our operators have constant coefficients, and hence commute. Namely, $e^{A+B/(1-q)}$, where $A$ and $B$ are regular at $q=1$, can be rewritten as $e^A e^{B/(1-q)}$.   

We are given that $f_{(\zeta)}(q^{1/m}/\zeta)=p\, \Psi^m (f_{(1)})$ where $p\in \hat{\K}_{+}$. Since $[q^{Q\p_Q},Q]=(q-1)Qq^{Q\p_Q}$ is divisible by $q-1$, for any finite-difference operator $B$, the commutator $\operatorname{ad}_B(p)=[B,p]$ with the operator of multiplication by $p$ is divisible by
$q-1$. Therefore $e^{B/(1-q)}p=Pe^{B/(1-q)}$, where $P=e^{\operatorname{ad}_{B/(1-q)}}(p)$ is regular at $q=1$. Thus, for some $P$ regular at $q=1$ we have:
\[ g_{(\zeta)}(q^{1/m}/\zeta) = e^{\Delta} P e^{\sum_{l>0}\la^{ml}\Psi^{ml}\left(D(1-q^{lQ\p_Q}, q^{\pm 1/m})\right)/ml(1-q^l)}\, \Psi^m (f_{(1)}).\]
  Comparing this expression with $\Psi^m(g_{(1)})$, take into account that $q^{\pm 1/m}$ coincides with $q^{\pm 1}$ modulo $q-1$, and $1/(1-q^{-lm})-1/m(1-q^{-l})$ is regular at $q=1$. Thus, again factoring off the terms regular at $q=1$, we conclude that  $g_{(\zeta)}(q^{1/m}/\zeta) $ is obtained from $\Psi^m(g_{(1)})$ by the application of an operator regular at $q=1$.  

From the explicit description of $\L^{fake}$, we have $g_{(1)}\in e^{\tau/(1-q)}\hat{\K}_{+}$ for some $\tau$. Therefore $\Psi^m(g_{(1)})\in e^{\Psi^m(\tau)/m(1-q)} \hat{\K}_{+}$. The latter is a $\D_q$-module, and hence  $g_{(\zeta)}(q^{1/m}/\zeta) \in \Psi^m(g_{(1)})\, \hat{\K}_{+}$ as required. 

Finally, for $\zeta \neq 0,\infty$, which is not a root of unity, regularity of $g$ at $q=\zeta^{-1}$ is obvious whenever the same is true for $f$. $\square$ 

\section*{$\Gamma$-operators}

{\tt Lemma.}{\em Let $l$ be a positive integer. Suppose that $\sum_{d\geq 0} f_dQ^d$ represents a point on the cone $\L \subset \K$. Then the same is true about:
\[ \sum_{d\geq 0} f_d\, Q^d\prod_{r=0}^{ld-1}(1-\la q^{-r}),\ \sum_{d\geq 0} \frac{f_d\, Q^d}{\prod_{r=1}^{ld}(1-\la q^r)}, \ \text{and}\  \sum_{d\geq 0}f_d\, Q^d\prod_{r=1}^{ld}(1-\la q^r). \]}

{\tt Proof.} We use $q$-Gamma-function  
\[ \Gamma_{q}(x) := e^{\sum_{k>0} x^k /k(1-q^k)} \sim \prod_{r=0}^{\infty} \frac{1}{1-xq^r}\]
for symbols of $q$-difference operators:
\begin{align*} \frac{\Gamma_{q^{-1}}(\la q^{-lQ\p_Q})}{\Gamma_{q^{-1}}(\la)}\, Q^d &= Q^d\, \frac{\prod_{r=-\infty}^0(1-\la q^r)}{\prod_{r=-\infty}^{-ld}(1-\la q^r)}  =Q^d\, \prod_{r=0}^{ld-1}(1-\la q^{-r}), \\
  \frac{\Gamma_{q^{-1}}(\la q^{lQ\p_Q})}{\Gamma_{q^{-1}}(\la)}\, Q^d &= Q^d\,  \frac{\prod_{r=-\infty}^0(1-\la q^r)}{\prod_{r=-\infty}^{ld}(1-\la q^r)} = \frac{Q^d} {\prod_{r=1}^{ld}(1-\la q^r)}, \ \text{and}\\
\frac{\Gamma_{q^{-1}}(\la)}{\Gamma_{q^{-1}}(\la q^{lQ\p_Q})}\, Q^d &= Q^d\,  \prod_{r=1}^{ld}(1-\la q^r) \ \text{respectively.} \end{align*}
The result follows now from the theorem of the previous section. $\square$

\section*{Application to fixed point localization}

In Part II, we used fixed point localization to characterize the range (denote it $\L_X$) of the big J-function in permutation- (and torus-) equivariant quantum K-theory of $X=\CC P^N$. Namely {\em a vector-valued ``rational function''  $f(q) = \sum_{i=0}^N f^{(i)}(q) \phi_i$ represents a point of $\L_X$ if and only if its components pass two tests, (i) and (ii):

  {\bf (i)} When expanded as meromorphic functions with poles $q\neq 0,\infty$  only at roots of unity, $f^{(i)} \in \L$, i.e. represent values of the big J-function $\J_{pt}$ in permutation-equivariant theory of the point target space;
  
 {\bf (ii)} Away from $q=0,\infty$, and roots of unity, $f^{i)}$ may have at most simple poles at $q = (\gL_j/\gL_i)^{1/m}$, $j\neq i$, $m=1,2,\dots$, with the residues  satisfying the recursion relations
  \[\Res_{q=(\gL_j/\gL_i)^{1/m}} f^{(i)}(q) \frac{dq}{q}=-\frac{Q^m}{C_{ij}(m)}\, f^{(j)}((\gL_j/\gL_i)^{1/m}),\]
where $C_{ij}(m)$ are explicitly described rational functions.}

We even verified that the hypergeometric series 
\[ J^{(i)}=(1-q) \sum_{d\geq 0}\frac{Q^d}{\left(\prod_{r=1}^d(1-q^r)\right) \prod_{j\neq i}\prod_{r=1}^d(1-q^r\gL_i/\gL_j)}\]
pass test (ii). Now we are ready for test (i). Indeed, we know
from Part I (or from Part III) that
\[ (1-q) \Gamma_q(Q):=(1-q)e^{\sum_{k>0}Q^k/k(1-q^k)}=(1-q)\sum_{d\geq 0}\frac{Q^d}
   {\prod_{r=1}^d(1-q^r)}\]
   lies in $\L$. According to Lemma, 
   \[ J^{(i)}=\prod_{j\neq i} \frac{\Gamma_{q^{-1}}(\gL_i\gL_j^{-1} q^{Q\p_Q})}{\Gamma_{q^{-1}}(\gL_j\gL_j^{-1})}\, (1-q)\Gamma_q(Q)\]
   also lies in $\L$. Thus, we obtain

   \medskip

   {\tt Corollary 1.} {\em The $K^0(\CC P^N)$-valued function
     \[ J_{\CC P^N} := \sum_{i=0}^N J^{(i)}\psi_i = (1-q)\sum_{d\geq 0}\frac{Q^d}{\prod_{j=0}^N\prod_{r=1}^d(1-P\gL_j^{-1}q^r)},\]
where $P=\O (-1)$ satisfies $\prod_{j=0}^N(1-P\gL_j^{-1})=0$,
represents a value of of the big J-function $\J_{\CC P^N}$.}

\medskip

{\tt Remark.} Note that all summands with $d>0$ are reduced rational functions of $q$, and so the Laurent polynomial part of $J_{\CC P^N}$ consists of the dilaton shift term $1-q$ only. This means that $J_{\CC P^N}$ represents the value of
   the big J-function $\J_{\CC P^N}(\t)$ at the input $\t=0$. Hence it is the small J-function (not only in permutation-equivariant but also) in the {\em ordinary} quantum K-theory of $\CC P^N$. In this capacity it was computed in \cite{GiL} by {\em ad hoc} methods. 
  
One can derive this way many other applications. To begin with, consider quantum K-theory on the target $E$ which is the total space of a vector bundle $E\to X$. To make the theory formally well-defined, one equips $E$ with the fiberwise scaling action of a circle, $T'$, and defines correlators by localization to fixed points $E^{T'}=X$ (the zero section of $E$).
This results in systematic {\em twisting} of virtual structure sheaves on the moduli spaces $X_{g,n,d}$ as follows:
\[ \O^{virt}_{g,n,d}(E):= \frac{\O^{virt}_{g,n,d}(X)}{\operatorname{Euler}^K_{T'}(E_{g,n,d})},\ \ E_{g,n,d}=(\ft_{n+1})_*\ev^*_{n+1}(E),\]
where the $T'$-equivariant K-theoretic Euler class of a bundle $V$ is defined by \[ \operatorname{Euler}^K_{T'}(V):=\tr_{\la \in T'} \left(\sum_k (-1)^k\bigwedge\!^kV^*\right) .\]
The division is possible in the sense that the $T'$-equivariant Euler class is invertible over the field of fractions of the group ring of $T'$. 
The elements $E_{g,n,d}\in K^0(X_{g,n,d})$ are invariant under permutations of the marked points. (In fact \cite{Co, CGi}, for $d\neq 0$,  $E_{g,n,d}=\ft^*E_{g,0,d}$
where $\ft: X_{g,n,d}\to X_{g,n,d}$ forgets all marked points.) Thus, we obtain a well-defined permutation-equivariant quantum K-theory of $E$.

\medskip

{\tt Corollary 2.} {\em Let $X=\CC P^N$, and $E=\oplus_{j=1}^M \O (-l_j)$.
Then the following $q$-hypergeometric series 
  \[ I_E := (1-q)\sum_{d\geq 0} \frac{Q^d}{\prod_{j=0}^N \prod_{r=1}^d(1-P\gL_j^{-1}q^r)} \prod_{j=1}^M\frac{\prod_{r=-\infty}^{l_jd-1}(1-\la P^{-l_j}q^{-r})}{\prod_{r=-\infty}^{-1}(1-\la P^{-l_j}q^{-r})}\]
represents a value of the big J-function in the permutation-equivariant quantum K-theory of $E$.}   

\medskip

Here $\la \in T'=\CC^{\times}$ acts on the fibers of $E$ as multiplication by $\la^{-1}$.  The K-theoretic Poincar{\'e} pairing on $X$ is twisted into $(a,b)_E=\chi (X; ab / \Eu^K_T(E))$.

\medskip

{\tt Example.} Let $X=\CC P^1$, $E=\O(-1)\oplus \O(-1)$. In $I_E$, pass to
the non-equivariant limit $\gL_0=\gL_1=1$:
\begin{align*} I_{E}=&(1-q)+ (1-\la P^{-1})^2 \times \\
  &(1-q)\sum_{d>0} Q^d \frac{(1-\la P^{-1}q^{-1})^2\cdots (1-\la P^{-1}q^{1-d})^2}{(1-Pq)^2(1-Pq^2)^2 \cdots (1-Pq^d)^2}.\end{align*}
The factor $(1-\la P^{-1})^2$, equal to $\Eu^K_{T'}$, reflects the fact that the part with $d>0$ is a push-forward from $\CC P^1$ to $E$. In the second non-equivariant limit, $\la=1$, it would turn into $0$ (since $(1-P^{-1})^2=0$ in
$K^0(\CC P^1)$. However, what the part with $d>0$ is push-forward of, survives in this limit:  
\[ (1-q)\sum_{d>0} \frac{Q^d}{P^{2d-2} q^{d(d-1)}(1-Pq^d)^2}, \ \text{where $(1-P)^2=0$.} \]
This example is usually used to extract information about ``local'' contributions of a rational curve $\CC P^{-1}$ lying in a Calabi-Yau 3-fold
with the normal bundle $\O(-1)\oplus \O(-1)$.

Note that decomposing the terms of this series into two summands: with poles at roots of unity, and with poles at $0$ or $\infty$, we obtain non-zero Laurent polynomials in each degree $d$. They form the input $\t=\sum_{d>0} \t_d(q,q^{-1})Q^d$ of the big J-function whose value $\J_E(\t)$ is given by the series. 

Finally, note that though the input is non-trivial, it is defined over the $\la$-algebra $\QQ [[Q]]$. This means that, although we are talking about permutation-equivariant quantum K-theory, the hypergeometric functions here, and in Corollary 2 in general, represent {\em symmetrized} K-theoretic GW-invariant, i.e. $S_n$-invariant part of the sheaf cohomology.

\medskip

Similarly, one can introduce K-theoretic GW-invariants of the {\em super-bundle}
$\Pi E$ (which is obtained from $E\to X$ by the ``parity change'' $\Pi$ of the fibers) by redefining the virtual structure sheaves as
\[ \O^{virt}_{g,n,d}(\Pi E) := \O^{virt}_{g,n,d}(X) \, \Eu^K_{T'}(E_{g,n,d}).\]
  When genus-$0$ correlators of this theory have non-equivariant limits (e.g. when $E$ is a positive line bundle, and $d>0$), the limits coincide with the appropriate correlators of the submanifold $Y\subset X$ given by a holomorpfic section of $\Pi E$.

\medskip

{\tt Corollary 3.} {\em Let $X=\CC P^N$, and $E=\oplus_{j =1}^M \O (l_j)$.
Then the following $q$-hypergeometric series 
  \[ I_{\Pi E} := (1-q)\sum_{d\geq 0} \frac{Q^d}{\prod_{j=0}^N \prod_{r=1}^d(1-P\gL_j^{-1}q^r)} \prod_{j=1}^M\frac{\prod_{r=-\infty}^{l_jd}(1-\la P^{l_j}q^r)}{\prod_{r=-\infty}^{0}(1-\la P^{l_j}q^r)}\]
represents a value of the big J-function in the permutation-equivariant quantum K-theory of $E$.}  

\medskip

Here $\la \in T'=\CC^{\times}$ acts on fibers of $E$ as multiplication by $\la$. The Poincar{\'e} pairing is twisted into $(a,b)_{\Pi E}=\chi (X; ab \Eu^K_T(E))$.

\medskip

{\tt Example.} When all $l_j>0$, it is safe pass to the non-equivariant limit $\gL_j=1$ and $\la=1$:
\[ I_{\Pi E}=(1-q)\sum_{q\geq 0} Q^d \frac{\prod_{j=1}^M\prod_{r=1}^{l_jd}(1- P^{l_j}q^r)}{\prod_{r=1}^d(1-Pq^r)^{N+1}},\]
which represents a value of the big J-function of $Y\subset \CC P^N$, pushed-forward from $K^0(Y)$ to $K^0(\CC P^N)$. Here $Y$ is a codimension-$M$ complete intersection given by equations of degrees $l_j$. Taking in account the degeneration of the Euler class in this limit, one may assume that $P$ satisfies the relation $(1-P)^{N+1-M}=0$.

When $\sum_j {l_j^2}\leq N+1$, the Laurent polynomial part of this series is
$1-q$, i.e. the corresponding input $\t$ of the J-function vanishes. In this case the series represents the small J-function of the ordinary quantum K-theory on $Y$. This result was obtained in \cite{GiT} in a different way:
based on the adelic characterization of the whole theory, but without the use of fixed point localization. As we have seen here, when $\t\neq 0$, the
series still represents the value $\J_Y(\t)$ in the {\em symmetrized} quantum K-theory of $Y$.

\medskip

In Part V these results will be carried over to all toric manifolds $X$, toric bundles $E\to X$, or toric super-bundles $\Pi E$. 
In fact, the intention to find a home for toric $q$-hypergeometric functions with non-zero Laurent polynomial part  was one of the motivations for developing the permutation-equivariant version of quantum K-theory.

\enddocument